\documentclass[11pt,reqno]{amsart}
\usepackage{tikz,amsmath,amsfonts,amscd,amssymb,epsf, euscript,mathtext}
\usepackage{amsxtra,mathrsfs,textcomp}
\usepackage{pdfsync}
\usepackage[all]{xy}
\usepackage{verbatim,epsf}
\newtheorem{theorem}{Theorem}

\theoremstyle{definition}

\theoremstyle{remark}

\newcommand{\field}[1]{\ensuremath{\mathbb{#1}}}
\newcommand{\CC}{\field{C}}

\newcommand{\HH}{\field{H}}
\newcommand{\PP}{\field{P}}
\newcommand{\RR}{\field{R}}

\newcommand{\ZZ}{\field{Z}}

\DeclareMathOperator{\im}{Im}
\newcommand{\del}{\partial}

\newcommand{\curly}[1]{\mathscr{#1}}

\newcommand{\cM}{\curly{M}}

\newcommand{\cP}{\curly{P}}
\newcommand{\cQ}{\curly{Q}}

\newcommand{\cS}{\curly{S}}

\newcommand{\NN}{\field{N}}

\usepackage{hyperref}
\hypersetup{colorlinks,citecolor=blue,plainpages=false,hypertexnames=false}
\begin{document}
\title[Real projective connections and the Liouville equation]{On real projective connections, V.I. Smirnov's approach, and black hole type solutions of the Liouville equation}
\author{Leon A. Takhtajan}
\address{Department of Mathematics,
Stony Brook University, Stony Brook, NY 11794-3651, USA;
The Euler International Mathematical Institute, Saint Petersburg, Russia}
\email{leontak@math.sunysb.edu}
\dedicatory{Dedicated to my teacher Ludwig Dmitrievich Faddeev on the occasion of his 80th birthday}
\maketitle
\begin{abstract}We consider real projective connections on Riemann surfaces and their corresponding solutions of the Liouville equation. We show that these solutions have singularities of special type (a black-hole type) on a finite number of simple analytical contours. We analyze the case of the Riemann sphere with four real punctures, considered in V.I. Smirnov's thesis (Petrograd, 1918), in detail.
\end{abstract}

\section{Introduction}
One of the central problems of mathematics in the second half of the 19th century and at the beginning of the 20th century was the problem of uniformization of Riemann surfaces. The classics, Klein  \cite {Klein} and Poincar\'{e} \cite{P1}, associated it with studying second-order ordinary differential equations with regular singular points. Poincar\'{e} proposed another approach to the uniformization problem \cite {P2}. It consists in finding a complete conformal metric of constant negative curvature, and it reduces to the global solvability of the Liouville equation, a special nonlinear partial differential equations of elliptic type on a Riemann surface.

Here, we illustrate the relation between these two approaches and describe solutions of the Liouville equation corresponding to second-order ordinary differential equations with a real monodromy group. In the modern physics literature on the Liouville equation it is rather commonly assumed that for the Fuchsian uniformization of a Riemann surface it suffices to have a second-order ordinary differential equation with a real monodromy group. But the classics already knew that this is not the case, and they analyzed  second-order ordinary differential equations with a real monodromy group on genus $0$ Riemann surfaces with punctures in detail. Nonetheless, they did not consider the relation to the Liouville equation, and we partially fill this gap here.

Namely, in Section \ref{proj}, following the lectures \cite{Turin}, we briefly describe the theory of projective connections on a Riemann surface --- an invariant method for defining a corresponding second-order ordinary differential equation with regular singular points. Following \cite{TZ1}, \cite{TZ2}, we review the main results on the Fuchsian uniformization, Liouville equation and the complex geometry of the moduli space. In Section \ref{real}, following \cite{G1}, we present the modern classification of projective connections with a real monodromy group, and review the results of V.I. Smirnov thesis \cite{S1} (Petrograd, 1918). This work, published in \cite{S2}, \cite{S3}, was the first where a complete classification of equations with a real monodromy group was given in the case of four real punctures. In Section \ref {real-2}, we give a modern interpretation of V.I. Smirnov results. Finally, in Section \ref{L-b-hole}, we describe solutions of the Liouville equation with black-hole type singularities associated with real projective connections. To the best of our knowledge, these solutions have not been considered previously.

\section{Projective connections, uniformization and the Liouville equation} \label{proj}
\subsection{Projective connections}  Let $X_{0}$ be a compact genus $g$ Riemann surface with marked points $x_{1},\dots,x_{n}$,  where $2g+n-2>0$, and let $\{U_{\alpha},z_{\alpha}\}$ be a complex-analytic atlas with local coordinates
$z_{\alpha}$ and transition functions  $z_{\alpha}=g_{\alpha\beta}(z_{\beta})$ on $U_{\alpha}\cap U_{\beta}$. Let
$X=X_{0}\setminus\{x_{1},\dots,x_{n}\}$ denote a corresponding Riemann surface of type $(g,n)$, a genus $g$ surface with $n$ punctures. The collection $R=\{r_{\alpha}\}$, where  $r_{\alpha}$ are holomorphic functions on $U_{\alpha}\cap X$, is called  \emph{(holomorphic) projective connection} on $X$, if on every intersection  $U_{\alpha}\cap U_{\beta}\cap X$,
$$r_{\beta}=r_{\alpha}\circ g_{\alpha\beta}(g'_{\alpha\beta})^{2}+\cS(g_{\alpha\beta}),$$
where $\cS(f)$ is the Schwarzian derivative of a holomorphic function $f$,
$$\cS(f)=\frac{f'''}{f'}-\frac{3}{2}\left(\frac{f''}{f'}\right)^{2}.$$
In addition, it is assumed that if $x_{i}\in U_{\alpha}$ and $z_{\alpha}(x_{i})=0$, then
\begin{equation} \label{order}
r_{\alpha}(z_{\alpha})=\frac{1}{2z_{\alpha}^{2}}+O\left(\frac{1}{|z_{\alpha}|}\right),\quad z_{\alpha}\rightarrow 0.
\end{equation} 
Projective connections form an affine space $\cP(X)$ over the vector space $\cQ(X)$ of holomorphic quadratic differentials on $X$;  elements of $\cQ(X)$ are collections $Q=\{q_{\alpha}\}$ with the transformation law 
$$q_{\beta}=q_{\alpha}\circ g_{\alpha\beta}(g'_{\alpha\beta})^{2}$$
and the additional condition that $q_{\alpha}(z_{\alpha})=O(|z_{\alpha}|^{-1})$ as $z_{\alpha}\rightarrow 0$, if $x_{i}\in U_{\alpha}$ and $z_{\alpha}(x_{i})=0$. The vector space $\cQ(X)$ has the complex dimension $3g-3+n$. (For more details on projective connections and quadratic differentials, see \cite{Turin} and the references therein).

A projective connection $R$ naturally determines a second-order linear differential equation on the Riemann surface $X$, the Fuchsian differential equation
\begin{equation} \label{Fuchs}
\frac{d^{2}u_{\alpha}}{dz_{\alpha}^{2}}+\frac{1}{2}r_{\alpha}u_{\alpha}=0,
\end{equation}
where $U=\{u_{\alpha}\}$ is understood and as a multi-valued  differential of order $-\frac{1}{2}$ on $X$. 
Equation \eqref{Fuchs} determines the monodromy group, a representation of the fundamental group $\pi_{1}(X,x_{0})$ of the Riemann surface $X$ with the marked point $x_{0}$ in $\mathrm{PSL}(2,\CC)=\mathrm{SL}(2,\CC)/\{\pm I\}$. Condition \eqref{order} implies that  the standard generators of $\pi_{1}(X,x_{0})$, which correspond  to the loops around the punctures $x_{i}$, are mapped to parabolic elements in $\mathrm{PSL}(2,\CC)$ under the monodromy representation.

\subsection{Uniformization} According to the uniformization theorem 
\begin{equation}\label{uniform}
X\cong\Gamma\backslash\HH,
\end{equation}
where $\HH=\{\tau\in\CC: \im\tau>0\}$ is the Poincar\'{e} model of the Lobachevsky plane, and $\Gamma\subset \mathrm{PSL}(2, \RR)$ is a type $(g,n)$ Fuchsian group acting on $\HH$ by fractional linear transformations. In other words, there exists a complex-analytic covering $J: \HH\rightarrow X$ whose automorphism group is $\Gamma$. The function inverse to $J$, a multi-valued analytic function $J^{-1}: X\rightarrow\HH$, is a locally univalent linear polymorphic function on $X$ (this means that its branches are connected by fractional linear transformations in $\Gamma$).
The Schwarzian derivatives of $J^{-1}$ with respect to $z_{\alpha}$ are well-defined on $U_{\alpha}$ and determine the \emph{Fuchsian projective connection} $R_{\mathrm{F}}=\{\cS_{z_{\alpha}}(J^{-1})\}$ 
on $X$, and the multi-valued functions $\displaystyle{\frac{1}{\sqrt{(J^{-1})'}}\;\text{and}\;\frac{J^{-1}}{\sqrt{(J^{-1})'}}}$ satisfy Fuchsian differential equation \eqref{Fuchs} with $R=R_{\mathrm{F}}$. The monodromy group of this equation, up to the conjugation in  $\mathrm{PSL}(2, \RR)$, is the Fuchsian group $\Gamma$.

Klein \cite{Klein} and  Poincar\'{e} \cite{P1} were solving the problem of uniformizing a Riemann surface $X$ by choosing a projective connection in Fuchsian equation \eqref{Fuchs} such that its monodromy group is a Fuchsian group $\Gamma$ with the property that \eqref{uniform} holds. But a direct proof of the existence of the Fuchsian projective connection $R_{\mathrm{F}}$  on $X$ turned out to be a very difficult problem, which has not been completely solved to this day (see \cite{H1}, \cite{H2}). In the case of type $(0,n)$ Riemann surfaces, to which we further restrict ourselves, this problem is formulated as follows. 

Let $X_{0}=\PP^{1}$ be the Riemann sphere and $X$ be a genus $0$ Riemann surface with $n$ punctures $z_{1},\dots,z_{n}$. Without loss of generality, we can assume that $z_{n-2}=0, z_{n-1}=1, z_{n}=\infty$ and $X=\CC\setminus\{z_{1},\dots,z_{n-3},0,1\}$. Equation \eqref{Fuchs} becomes
\begin{equation} \label{genus-0}
\frac{d^{2}u}{dz^{2}}+\frac{1}{2}\sum_{i=1}^{n-1}\left(\frac{1}{2(z-z_{i})^{2}}+\frac{c_{i}}{z-z_{i}}\right)u=0,
\end{equation}
where $z$ is a global complex coordinate on $X$.  The complex parameters $c_{1},\dots,c_{n-1}$ satisfy the two conditions
\begin{equation} \label{acc}
\sum_{i=1}^{n-1}c_{i}=0,\quad\sum_{i=1}^{n-1}z_{i}c_{i}=1-\frac{n}{2},
\end{equation}
which allow expressing $c_{n-2}$ and $c_{n-1}$ explicitly in terms of  $z_{1},\dots,z_{n-3}$ 
and the remaining $n-3$ parameters $c_{1},\dots,c_{n-3}$.

In the classical approach of Klein and Poincar\'{e} to the uniformization
problem, given the singular points $z_{1},\dots,z_{n-3},0,1,\infty$, it was required to
choose parameters $c_{1},\dots, c_{n-3}$ such that the
monodromy group of equation \eqref{genus-0} is a Fuchsian group isomorphic to the fundamental group of the Riemann surface $X$.
The ratio of two linear independent solutions of equation \eqref{genus-0} up to a fractional linear transformation would then be the desired multi-valued mapping $J^{-1}: X\rightarrow\HH$ realizing the uniformization of the Riemann surface $X$. The corresponding $\Gamma$-automorphic function $J:\HH\rightarrow\CC$ is called \emph{Klein's Hauptmodul} (Hauptfunktion). The complex numbers $c_{1},\dots, c_{n-3}$, the \emph{accessory parameters} of the Fuchsian uniformization of the surface $X$, are uniquely determined by the singular points  $z_{1},\dots,z_{n-3}$. Moreover
\begin{equation} \label{J-inverse}
\cS(J^{-1})(z)=\sum_{i=1}^{n-1}\left(\frac{1}{2(z-z_{i})^{2}}+\frac{c_{i}}{z-z_{i}}\right).
\end{equation}

To prove the existence of the accessory parameters, Poincar\'{e} proposed the so-called continuity method  in \cite{P1}. But a rigorous solution of the uniformization problem could not be obtained using this method, and being subjected to criticism, the method was soon abandoned.
The ultimate solution of the uniformization problem was obtained by Koebe and Poincar\'{e} in 1907 by using quite different methods, in particular, by using potential theory (see, e.g.,~\cite{FK} for a modern exposition).

 \subsection{The Liouville equation} \label{Liouv}
The projection on $X$ of the Poincar\'{e} metric $(\im\tau)^{-2}|d\tau|^{2}$ on $\HH$ is a complete conformal metric on $X$ of constant negative curvature $-1$. It has the form $e^{\varphi(z)}|dz|^{2}$, where
\begin{equation} \label{L-formula}
e^{\varphi(z)}=\frac{|(J^{-1})'(z)|^{2}}{(\im J^{-1}(z))^{2}}.
\end{equation}
The smooth function $\varphi$ on $X$ satisfies the Liouville equation
\begin{equation} \label{L}
\varphi_{z\bar{z}}=\frac{1}{2}e^{\varphi}
\end{equation}
and has the asymptotic behavior
\begin{equation} \label{as}
\varphi(z)=\begin{cases}-2\log|z-z_{i}|-2\log|\log|z-z_{i}||+ o(1),\;\;z\rightarrow z_{i},\; i\neq n,\\
-2\log|z|-2\log\log|z| +o(1),\;\;z\rightarrow\infty.
\end{cases}
\end{equation}
In \cite{P2}, Poincar\'{e} proposed an approach to the uniformization problem based on the Liouville equation. Namely, he proved that Liouville equation \eqref{L} is uniquely solvable in the class of smooth real-valued functions on $X$ with asymptotic behavior \eqref{as}. It hence follows that $T_{\varphi}=\varphi_{zz}-\frac{1}{2}\varphi_{z}^{2}$ is a rational function of the form \eqref{J-inverse} and that the differential equation 
$$\frac{d^{2}u}{dz^{2}}+\frac{1}{2}T_{\varphi}u=0$$
has a Fuchsian monodromy group that uniformizes the Riemann surface $X$.

The Liouville equation is the Euler-Lagrange equation for the functional
$$S(\psi)=\lim_{\varepsilon \rightarrow0}\left(\iint\limits_{X_{\varepsilon}}(|\psi_{z}|^{2}+e^{\psi})d^{2}z+2\pi n\log\varepsilon + 4\pi(n-2)\log|\log\varepsilon|\right),$$
where $d^{2}z$ is the Lebesgue measure on $\CC$, 
$$\displaystyle{X_{\varepsilon}=X\setminus\left(\bigcup_{i=1}^{n-1}\{|z-z_{i}|<\varepsilon\}\bigcup\{|z|>1/\varepsilon\}\right)},$$ 
and $\psi$ belongs to the class of smooth functions $X$ with asymptotic behavior \eqref{as}.
The quantity $T_{\psi}=\psi_{zz}-\frac{1}{2}\psi_{z}^{2}$ plays the role of the $(2,0)$-component of the stress-energy tensor in the classical Liouville theory, and 
$$T_{\varphi}=\cS(J^{-1}).$$

We let
$$\cM_{0,n}=\{(z_{1},\dots,z_{n-3})\in\CC^{n-3}: z_{i}\neq z_{j}\;\;\text{for}\;\; i\neq j,\;\;\text{and}\;\;
z_{i}\neq 0,1\}$$
denote the moduli space of genus $0$ Riemann surfaces with $n$ ordered punctures (rational curves with $n$ marked points).
The critical values of the Liouville action functional $S(\psi)$ (the values on the extrema $\varphi$ for all surfaces $X$) determine a smooth function $S:\cM_{0,n}\rightarrow\RR$, the \emph{classical action for the Liouville equation}. As proved in \cite{TZ1}, \cite{TZ2},  the classical action for the Liouville equation plays a fundamental role in the complex geometry of the moduli space $\cM_{0,n}$. Namely,
$S$ is a common antiderivative for the accessory parameters
$$c_{i}=-\frac{1}{2\pi}\frac{\del S}{\del z_{i}},\quad i=1,\dots,n-3,$$
and also a K\"{a}hler potential for the Weil-Petersson metric on $\cM_{0,n}$,
$$-\frac{\del^{2}S}{\del z_{i}\del\bar{z}_{j}}=\left\langle\frac{\del}{\del z_{i}},\frac{\del}{\del z_{j}}\right\rangle_{WP},\quad i,j=1,\dots n-3.$$

The statement that the classical action for the Liouville equation is a common antiderivative for accessory parameters for genus $0$ surfaces was conjectured by Polyakov\footnote{Lecture at Leningrad's branch of V.A. Steklov Mathematical Institute, 1982, unpublished.} based of the semiclassical analysis of the conformal Ward identities of the quantum Liouville theory (see \cite{T1}). 

\section{Real projective connections and V.I. Smirnov thesis} \label{real}
\subsection{General case} 
Let $X$ be a Riemann surface of type $(g,n)$. A projective connection $R$ on $X$ is said to be  \emph{real} or \emph{Fuchsian} if its monodromy group up to a conjugation in $\mathrm{PSL}(2,\CC)$ is respectively a subgroup in $\mathrm{PSL}(2,\RR)$ or a Fuchsian group. By the uniformization theorem, a Fuchsian projective connection $R_{\mathrm{F}}$ is uniquely characterized by the condition that its monodromy group is precisely the Fuchsian group $\Gamma$ that uniformizes the Riemann surface $X$. It is rather natural to ask whether it is possible to characterize a projective connection $R_{\mathrm{F}}$ on $X$ by simpler conditions like being real (see \cite[p. 224]{F1}) or Fuchsian? The answer to the question is negative.

Namely, for a compact genus $g>1$ Riemann surface  $X=\Gamma\backslash\HH$,  Goldman \cite{G1} showed that to every
integral Thurston's measurable geodesic lamination $\mu=\sum_{i}m_{i}\gamma_{i}$ (disjoint union), where $\gamma_{i}$ are simple closed geodesics in the hyperbolic metric on $X$, and $m_{i}$ are non-negative integers, there corresponds a  genus $g$ Riemann surface  $Gr_{\mu}(X)$  with a projective connection $R(\mu)$ having the monodromy group  $\Gamma$. Riemann surfaces $Gr_{\mu}(X)$ are obtained from $X$ by the so-called ``grafting'' procedure that generalizes classic examples of Maskit-Hejhal and Sullivan-Thurston (see \cite{G1}). Moreover, the set of all Fuchsian projective connections on all genus $g$ Riemann surfaces  is isomorphic to the direct product of the Teichm\"{u}ller space $T_{g}$ and the set of integral Thurston's measurable geodesic laminations on a genus $g$ topological surface. It was proved in \cite{Tan1} that on each genus $g>1$ Riemann surface $X$ here are infinitely many Fuchsian projective connections.

Real projective connections on Riemann surfaces of type $(g, n)$ were also studied by Faltings \cite{F1}. As shown for compact Riemann surfaces in \cite{G1}, to each half-integral Thurston's measurable lamination $\mu$ there is a genus $g$ Riemann surface $Gr_ {\mu}(X)$ with a real projective connection.

\subsection{Surfaces of type $(0,4)$ and V.I. Smirnov's approach} \label{real-2}
The classics associated the uniformization problem of Riemann surfaces with differential equations. As a basic example, they considered the case of Riemann surfaces of type $(0,4)$; the corresponding problem was to find an accessory parameter in equation \eqref {genus-0} such that its monodromy group was Fuchsian or Kleinian group (a discrete subgroup of $\mathrm{ PSL}(2, \CC)$). For the special case of real singular points, Klein \cite{K2} proposed an approach that uses Sturm's oscillation theorem, and Hilb proved \cite{Hilb} that  equation \eqref{genus-0} has a Fuchsian monodromy group for infinitely many values of the accessory parameter. Hilbert \cite[Kap. XX] {Hilbert} reduced this problem to the study of a certain integral equation.

The problem of a real monodromy group of equation \eqref{genus-0} with four real singular points was completely solved by V.I. Smirnov in his thesis \cite{S1}, published in Petrograd in 1918 (its main content was presented in \cite{S2}, \cite{S3}). Namely, we consider equation
 \eqref{genus-0} with the singular points $z_{1}=0, z_{2}=a, z_{3}=0$ and $z_{4}=\infty$, where $0<a<1$. Writing the general solution of equations  
\eqref{acc} in the form
$$c_{1}=1+\frac{1+2\lambda}{a},\;c_{2}=\frac{1+2\lambda}{a(a-1)},\;c_{3}=-\frac{a+2\lambda}{a-1},$$
where $\lambda$ is the accessory parameter, and changing the dependent variable $y=\sqrt{z(z-a)(z-1)}\,u$, we obtain the equation
\begin{equation} \label{Smirnov}
\frac{d}{dz}\left(p(z)\frac{dy}{dz}\right)+(z+\lambda)y=0,\quad p(z)=z(z-a)(z-1).
\end{equation}

Let $(y^{(1)}_{i},y^{(2)}_{i})$ denote the standard basis in the solution space of \eqref{Smirnov},  which in the neighborhood a singular point   $z_{i}$ consists of normalized holomorphic solutions
\begin{align*}
y^{(1)}_{i}(z,\lambda) & =1+\sum_{k=1}^{\infty}a_{ik}(z-z_{i})^{k},\quad i=1,2,3,\\
y^{(1)}_{4}(z,\lambda) & =\frac{1}{z}+\sum_{k=2}^{\infty}\frac{a_{4k}}{z^{k}} \\
\intertext{and}
y^{(2)}_{i}(z,\lambda)&=y^{(1)}_{i}(z,\lambda)\log(z-z_{i}) + \tilde{y}_{i}(z,\lambda),\quad i=1,2,3, \\
y^{(2)}_{4}(z,\lambda) &=y^{(1)}_{4}(z,\lambda)\log\frac{1}{z}+\tilde{y}_{4}(z,\lambda),
\end{align*}
where $\tilde{y}_{i}(z,\lambda)$ are holomorphic in the neighborhood of $z_{i}$.  For real $\lambda$ the power series $y^{(1,2)}_{i}(z,\lambda)$ and $\tilde{y}_{i}(z,\lambda)$ have real coefficients.

To determine a real $\lambda$ for which the monodromy group of equation  \eqref{Smirnov} is real, the classics used  the notion of  \emph{real continuation}.  Namely (see \cite{S1}), if we have  
$$y(z)=c\log(a-z) + f(z)$$
in the neighborhood of a singular point  $a$ for real  $z<a$, 
where $\log 1=0$ and the function $f(z)$ is holomorphic in a neighborhood of $a$, then the real continuation of $y$ to the domain $z>a$ is defined as
$$y(z)=c\log(z-a)+f(z).$$
The following statement holds:
\begin{theorem}[Klein, Hilbert] Equation \eqref{Smirnov} has a real monodromy group if $\lambda$ is real and one of the following conditions holds.
\begin{itemize}
\item[1.] The solution $y^{(1)}_{0}(z,\lambda)$ is holomorphic in a neighborhood of the singular point $z_{2}=a$. 
\item[2.] The solution  $y^{(1)}_{2}(z,\lambda)$ is holomorphic in a neighborhood of the singular point $z_{3}=1$.
\item[3.] Under the real continuation through $z=a$ the solution $y^{(1)}_{0}(z,\lambda)$ is holomorphic in a neighborhood of the singular point $z_{3}=1$.
\end{itemize}
\end{theorem} 
\noindent
Moreover, under conditions 1,2, and 3, the respective ratios $\displaystyle{\eta=\sqrt{-1}\,y_{3}^{(1)}/y_{1}^{(1)}}$, $\displaystyle{\eta=\sqrt{-1}\,y_{2}^{(1)}/y_{4}^{(1)}}$ and $\displaystyle{\eta=\sqrt{-1}\,y_{1}^{(1)}/y_{2}^{(1)}}$ of linear independent solutions of equation \eqref{Smirnov} are transformed by real fractional linear  transformations when going around the singular points $0,a$ and $1$ (see \cite{S1} for details).

Conditions 1-3 determine the following three Sturm-Liouville type spectral problems for equation  \eqref{Smirnov}. Namely, we must determine the values of $\lambda$ such that
\begin{itemize}
\item[1.] there is a solution on the interval $[0,a]$ that is regular at $0$ and $a$;
\item[2.] there is a solution on the interval $[a,1]$ that is regular at $a$ and $1$;
\item[3.] there is a solution regular at $0$ such that it is regular at $1$ under the real continuation through $a$. 
\end{itemize}  
Using the classical Sturm method, V.I. Smirnov \cite{S1} proved the following result.
\begin{theorem}[V.I. Smirnov] \label{VIS} Each of Sturm-Liouville problems  $1-3$ has a simple unbounded discrete spectrum. 
Namely, the following statements hold:
\begin{itemize}
\item[1.]  Spectral problem $1$ has infinitely many eigenvalues  $\mu_{k}$, $k\in\NN$, accumulating at $\infty$ and satisfying the inequalities
$$-a<\mu_{1}<\mu_{2}<\dots$$
\item[ii)] Spectral problem $2$ has infinitely many eigenvalues  $\mu_{-k}$, $k\in\NN$, accumulating at $-\infty$ and satisfying the inequalities
$$-a>\mu_{-1}>\mu_{-2}>\dots$$
\item[iii)] Spectral problem $3$ has infinitely many eigenvalues $\lambda_{k}$, $k\in\ZZ$, accumulating at $\pm\infty$ and satisfying the inequalities

$$\dots<\mu_{-2}<\lambda_{-1}<\mu_{-1}<\lambda_{0}<\mu_{1}<\lambda_{1}<\mu_{2}<\dots$$
\end{itemize} 
\end{theorem}

The case $\lambda=\lambda_{0}$ corresponds to the Fuchsian uniformization of the Riemann surface $X=\CC\setminus\{0,a,1\}$ and the ratio $\eta=\sqrt{-1}\,y_{1}^{(1)}/y_{2}^{(1)}$ bijectively maps the upper half-plane of $z$ to the interior of a circular rectangle with zero angles and sides orthogonal to $\RR\cup\{\infty\}$. Normalizing $\eta$ by a real fractional linear transformation such that the images of all singular points $0,a,1$ and $\infty$ are finite, we obtain the rectangle in Fig. 1  (\emph{cf}. Fig. 9 in \cite{K2}).

\vspace{4mm}
\begin{tikzpicture}[thick]
\draw[fill=gray!40!white] (0,0) to [out=90,in=180] (.7,0.7) to [out=0,in=90] (1.4,0) to [out=90,in=180] (2.3,.9) to [out=0,in=90] (3.2,0) to [out=90,in=180]
(3.7,.5) to [out=0,in=90] (4.2,0) to [out=90,in=0] (2.1,2.1) to [out=180,in=90] (0,0);
\draw (-1,0) -- (5,0);
\node at (0,0) {$\bullet$};
\node at (1.4,0) {$\bullet$};
\node at (3.2,0) {$\bullet$};
\node at (4.2,0) {$\bullet$};
\node at (.7,0.7) {$>$};
\node at (2.3,0.9) {$>$};
\node at (3.7,0.5) {$>$};
\node at (2.1,2.1) {$<$};
\node at (0,-.3) {$\eta(\infty)$};
\node at (1.4,-0.3){$\eta(0)$};
\node at (3.2,-.3) {$\eta(a)$};
\node at (4.2,-.3) {$\eta(1)$};
\node at (2.5,-1) {$\text{Fig. 1}$};
\end{tikzpicture}

\noindent
Analytically continuing $\eta(z)$ to the lower half-plane of $z$, we obtain a multi-valued linear polymorphic function $\eta: X\rightarrow\HH$  with a Fuchsian group $\Gamma$ such that $J=\eta^{-1}$ determines isomorphism \eqref{uniform}. 

The corresponding ratio $\eta$ is a bijective function on the upper half-plane of  $z$ also in the cases $\lambda=\mu_{\pm 1}$. Hence, if $\lambda=\mu_{1}$, then we have $\eta(0)=\eta(a)=\infty$ and $\eta(1)=\eta(\infty)=0$. Normalizing $\eta=\sqrt{-1}\,y_{3}^{(1)}/y_{1}^{(1)}$ such that the images of the singular points are finite, we obtain a bijective map $\eta$ of the upper half-plane of $z$ onto the interior of the degenerate circular rectangle in Fig. 2. The corresponding monodromy groups are  Schottky groups.

\vspace{4mm}
\begin{tikzpicture}[thick]
\draw[fill=gray!40!white] (0,0) to [out=90,in=180] (.7,0.7) to [out=0,in=90] (1.4,0) to [out=90,in=180] (2.3,.9) to [out=0,in=90] (3.2,0) to 
[out=-90,in=0] (2.3,-.9) to [out=180,in=-90] (1.4,0) to [out=-90,in=0] (.7,-.7) to [out=180,in=-90] (0,0) to [out=-90,in=180] (2.1,-2.1) to [out=0,in=-90] 
(4.2,0)  to [out=90,in=0] (2.1,2.1) to [out=180,in=90] (0,0);
\draw (-1,0) -- (5,0);
\node at (0,0) {$\bullet$};
\node at (1.4,0) {$\bullet$};
\node at (.7,0.7) {$>$};
\node at (2.3,0.9) {$>$};
\node at (2.1,2.1) {$<$};
\node at (2.4,-2.8) {$\text{Fig. 2}$};
\end{tikzpicture}

\noindent
For all other values of $\lambda_{k}$ and $\mu_{k}$ the corresponding map $\eta$ is no longer a bijective map on the upper half-plane of $z$. Hence, if $\lambda=\lambda_{1}$, then the upper half-plane is mapped onto the interior of the annulus in Fig. 3  (\emph{cf}. Fig. 10 in \cite{K2}). Here, the function $\eta$ takes the values twice from the marked darker domain, which corresponds to the rectangle in Fig. 1. When $\lambda=\lambda_{k}$, then this rectangle is wrapped over itself $2|k|$ times. 

\vspace{4mm}
\begin{tikzpicture}[thick]
\draw[fill=gray!80!white] (0,0) to [out=90,in=180] (.7,0.7) to [out=0,in=90] (1.4,0) to [out=90,in=180] (2.3,.9) to [out=0,in=90] (3.2,0) to [out=90,in=180]
(3.7,.5) to [out=0,in=90] (4.2,0) to [out=90,in=0] (2.1,2.1) to [out=180,in=90] (0,0);
\draw[fill=gray!40!white] (0,0) to [out=90,in=180] (.7,0.7) to [out=0,in=90] (1.4,0) to [out=-90,in=180] (2.3,.-.9) to [out=0,in=-90] (3.2,0) to [out=90,in=180]
(3.7,.5) to [out=0,in=90] (4.2,0) to [out=-90,in=0] (2.1,-2.1) to [out=180,in=-90] (0,0);
\draw (-1,0) -- (5,0);
\node at (0,0) {$\bullet$};
\node at (1.4,0) {$\bullet$};
\node at (3.2,0) {$\bullet$};
\node at (4.2,0) {$\bullet$};

\node at (2.4,-2.8) {$\text{Fig. 3}$};
\end{tikzpicture}

\noindent
Similarly, if $\lambda=\mu_{2}$, then the upper half-plane of  $z$ maps onto the interior of the annulus in Fig. 4. Here, the function 
$\eta$ takes the values twice from the marked darker domain, which corresponds to the degenerate rectangle in Fig. 2. If $\lambda=\mu_{k}$, then this rectangle is wrapped over itself $|k|$ times.

\vspace{4mm}
\begin{tikzpicture}[thick]
\draw[fill=gray!80!white] (0,0) to [out=90,in=180] (.7,0.7) to [out=0,in=90] (1.4,0) to [out=90,in=180] (2.3,.9) to [out=0,in=90] (3.2,0) to 
[out=-90,in=0] (2.3,-.9) to [out=180,in=-90] (1.4,0) to [out=-90,in=0] (.7,-.7) to [out=180,in=-90] (0,0) to [out=-90,in=180] (2.1,-2.1) to [out=0,in=-90] 
(4.2,0)  to [out=90,in=0] (2.1,2.1) to [out=180,in=90] (0,0);
\draw[fill=gray!40!white] (0,0) to [out=90,in=180] (.7,0.7) to [out=0,in=90] (1.4,0) to [out=-90,in=0] (.7,-.7) to [out=180,in=-90] (0,0);
\draw (-1,0) -- (5,0);
\node at (0,0) {$\bullet$};
\node at (1.4,0) {$\bullet$};
\node at (2.4,-2.8) {$\text{Fig. 4}$};
\end{tikzpicture}

\noindent
It is instructive to compare these results of V.I. Smirnov with Goldman's classification of Fuchsian and real projective connections on Riemann surfaces generalized for the surfaces of type $(g,n)$. The Fuchsian series $\lambda=\lambda_{k}$ corresponds to integral laminations in \cite{G1}, while series $\lambda=\mu_{k}$ corresponds to half-integral laminations.

\section{Black hole type solutions of the Liouville equation} \label{L-b-hole}

The Fuchsian uniformization of a Riemann surface $X$ determines a solution of the Liouville equation: a smooth function $\varphi$ on $X$, satisfying equation \eqref{L} and having asymptotic behavior \eqref{as} (see Section \ref{Liouv}). The function $\varphi$ is obtained from the ratio   $J^{-1}$ of linear independent solutions of equation \eqref{genus-0} by formula \eqref{L-formula}. This formula is well-defined because of the condition that the monodromy group $\Gamma$ of equation \eqref{genus-0} is real; the smoothness of $\varphi$ is ensured by the condition that $\Gamma$ uniformizes the Riemann surface $X$ and its image under the multi-valued map $J^{-1}$ is the upper half-plane $\HH$.  

Similarly, with each equation \eqref{genus-0}  with the real monodromy group one associates a solution of the Liouville equation. We set
\begin{equation} \label{L-general}
e^{\varphi(z)}=\frac{|\eta'(z)|^{2}}{\left(\mathrm{Im}\,\eta(z)\right)^{2}},
\end{equation}
where $\eta$ is the ratio of linear independent solutions of equation \eqref{genus-0}, which transforms by linear fractional transformations when going around the singular points ($\eta=J^{-1}$ in the Fuchsian case). The function $\varphi$ is well-defined because of the realness of the monodromy group and has asymptotic behavior \eqref{as}. The latter follows from the theory of Fuchsian equations with equal exponents. But solution \eqref{L-general} is no longer smooth: the image of $X$ under the multi-valued map $\eta$ has a nontrivial intersection with the real axis, and the function $\varphi$ is singular on $\eta^{-1}(\RR)$.

Namely, it follows from results in \cite[\S6]{F1} that the inverse image $\eta^{-1}(\RR)$ is a disjoint union of finitely many simple closed analytic curves on $X$. Let $C$ be one such curve. There is a branch of the multi-valued function $\eta$ that maps $C$ bijectively onto the circle, so that $C=\{z=\eta^{-1}(t), t\in[\alpha,\beta]\}$. It is convenient to introduce the  \emph{Schwarz function} $S$ of the analytic contour $C$ by the formula 
$$S=\bar\eta^{-1}\circ\eta,$$
where $\bar{\eta}^{-1}(z)=\overline{\eta^{-1}(\bar{z})}$. The Schwarz function is defined in some neighborhood of the contour  $C$ and determines it by the equation $\bar{z}=S(z)$ (see \cite{Davis}).  It is easy to show that in terms of the Schwarz function, the solution $\varphi$ has the same singularities on $C$ as the function
\begin{equation} \label{Schwarz-1}
-\frac{4\overline{S'(z)}}{(z-\overline{S(z)})^{2}}.
\end{equation}
Namely, as  $z\rightarrow z_{0}\in C$ along any direction non tangent to $C$,
\begin{equation} \label{Schwarz-2}
e^{\varphi(z)}=-\frac{4S'(z_{0})}{\left(\bar{z}-\bar{z}_{0}-S'(z_{0})(z-z_{0})\right)^{2}}(1+O(|z-z_{0}|)).
\end{equation}
Note that due to the condition $\bar{S}(S(z))=z$ the function in the right hand side of  \eqref{Schwarz-2} is real and positive. The singularities of type \eqref{Schwarz-1}--\eqref{Schwarz-2} on a contour $C$ are similar to the singularity on  $\RR$ of the Poincar\'{e} metric on $\CC\setminus\RR$, which corresponds to the Schwarz function $S(z)=z$.

We can therefore state the following problem: on the Riemann surface  $X=\CC\setminus\{z_{1},\dots,z_{n-3},0,1\}$, find simple analytic contours $C_{1},\dots,C_{k}$ and a function $\varphi$ such that on $X\setminus\cup_{j=1}^{k}C_{j}$ the function  $\varphi$ satisfy Liouville equation  \eqref{L}, has asymptotic behavior \eqref{as} at the punctures $z_{i}$ and has singularities of the type \eqref{Schwarz-1}, \eqref{Schwarz-2} on the contours $C_{j}$. On each connected components of $X\setminus\cup_{j=1}^{k}C_{j}$, $e^{\varphi(z)}|dz|^{2}$ determines a complete metric of constant negative curvature $-1$.  The boundary $C_{j}$ may be interpreted as the horizon of a black hole, and we therefore call corresponding solutions of the Liouville equation solutions of the \emph{black-hole type}. It follows from Goldman's classification of real projective connections \cite{G1} that there exists a family of such solutions parameterized by the ``integral lattice'' of integral and half-integral measurable Thurston's laminations, implicitly defined by the grafting procedure. 

From the results in V.I. Smirnov thesis, we obtain a rather explicit description of black-hole type solutions in the case of four real singular points. 
Namely,  we obtain the following result from Theorem \ref{VIS}.
 \begin{theorem} All black hole type solutions of the Liouville equation with four real punctures $0,a,1$ and $\infty$ are described as follows:
 \begin{itemize}
 \item[1)] Solutions of the Fuchsian type, which correspond to the values of the accessory parameter $\lambda=\lambda_{k}$ with integer $k$ and have $2|k|$ contours  $C_{j}$:  these contours go over the points $0$ and $a$ if $k>0$  and over the points $a$ and $1$ if $k<0$.
 \item[2)]  Solutions of the Schottky type, which correspond to the values of the accessory parameter $\lambda=\mu_{k}$ with integer $k\neq 0$ and have  $2|k|-1$ contours $C_{j}$:  these contours go over the points $0$ and $a$ if $k>0$ and over the points $a$ and $1$ if $k<0$.
 \end{itemize}
 \end{theorem} 
 
 In the general case, it is convenient to substitute $\chi(z)=e^{-\varphi(z)/2}$, which transform Liouville equation \eqref{L} into
\begin{equation}\label{L-chi}
-\chi\,\chi_{z\bar{z}}+|\chi_{z}|^{2}=\frac{1}{4},
\end{equation}
and asymptotics \eqref{as} into
\begin{equation} \label{as-chi}
\chi(z)=\begin{cases}|z-z_{i}||\log|z-z_{i}||(1+ o(1)),\;\;z\rightarrow z_{i},\; i\neq n,\\
|z|\log|z|(1 +o(1)),\;\;z\rightarrow\infty.
\end{cases}
\end{equation}
Singularities \eqref{Schwarz-1} transform into the vanishing condition on the contour $C$:
\begin{equation} \label{Schwarz-3}
\chi(z)\sim\frac{z-\overline{S(z)}}{2\sqrt{-\overline{}S'(z)}},
\end{equation} 
and the real-valued function $\chi(z)$ hence changes sign under the Schwarz reflection $z^{*}=S(z)$ through $C$. Elliptic partial differential equation  \eqref{L-chi} with asymptotic behavior \eqref{as-chi} and vanishing conditions  \eqref{Schwarz-3} on the contours $C_{j}$ is a boundary value problem with a free boundary.  It would interesting to use the method of continuation with respect to a parameter together with the a priori estimates to solve it, as it was done in \cite{P2} for Liouville equation \eqref{L} with asymptotic behavior \eqref{as}.

In conclusion, we note that the function $\chi$ plays an important role in the theory of the Liouville equation. Namely, it is a bilinear form in solutions of equation  \eqref{genus-0} and their complex conjugates and satisfies equation \eqref{genus-0} 
\begin{equation} \label{2}
\chi_{zz}+\frac{1}{2}T_{\varphi}\chi=0
\end{equation}
and the complex conjugate equation. In the quantum Liouville theory, the field  $\chi=e^{-\varphi(z)/2}$ describes a vector degenerate at the level $2$ in a Verma module for the Virasoro algebra. For the black hole type solutions the function $\chi$ still satisfies equation \eqref{2}. 
It would be interesting to elucidate what role it plays in the quantum Liouville theory.

\subsection *{Acknowledgements} The author is pleased to express his gratitude to P.G. Zograf for the useful discussion of the results in \cite{G1}, \cite{F1} and for the geometric interpretation of the results in V.I. Smirnov's thesis. The results in this paper were presented in part at the conferences `` Perspectives, Open Problems \& Applications of Quantum Liouville Theory'' at Stony Brook, USA in March 2010 and `` Mathematics - XXI century. 70 years PDMI'' in St. Petersburg in September 2010. This work was supported by the National Science Foundation (NSF, grants  DMS-0705263 and DMS-1005769).

 \end{document}